\documentclass[12pt,reqno]{amsart}

\usepackage{setspace}

\usepackage{caption}

\usepackage
{hyperref, epigraph}

\usepackage{breakurl}   

\usepackage{graphicx} 

\usepackage[export]{adjustbox}  

\usepackage{enumitem} 
\theoremstyle{definition}

\DeclareMathOperator{\R}{{\mathbb R}}

\numberwithin{equation}{section}

\author[Bair]{Jacques Bair}\address{J. Bair, HEC-ULG, University of
Liege, 4000 Belgium}\email{j.bair@ulg.ac.be}

\author[B\l{}aszczyk]{Piotr B\l{}aszczyk}\address{P. B\l{}aszczyk,
Institute of Mathematics, Pedagogical University of Cracow,
Poland}\email{pb@up.krakow.pl}

\author[Fuentes Guill\'en]{ El\'{\i}as Fuentes Guill\'en}
\address{E. Fuentes Guill\'en, Department of Mathematics, Faculty of
Sciences, UNAM, Mexico}\email{eliasfuentesguillen@gmail.com}

\author[Heinig]{Peter Heinig} \address{P. Heinig}
\email{peter.c.heinig@gmail.com}

\author[Kanovei]{Vladimir Kanovei} \address{V. Kanovei, IPPI RAS,
Moscow, Russia}\email{kanovei@googlemail.com}

\author[Katz]{Mikhail G. Katz}\address{M. Katz, Department of
Mathematics, Bar Ilan University, Ramat Gan 5290002
Israel}\email{katzmik@macs.biu.ac.il}

\begin{document}


\thispagestyle{empty}


\title[Continuity between Cauchy and Bolzano] {Continuity between
Cauchy and Bolzano: Issues of antecedents and priority}

\begin{abstract}

In a paper published in 1970, Grattan-Guinness argued that Cauchy, in
his 1821 book \emph{Cours d'Analyse}, may have plagiarized Bolzano's
book \emph{Rein analytischer Beweis} (RB), first published in 1817.
That paper was subsequently discredited in several works, but some of
its assumptions still prevail today. In particular, it is usually
considered that Cauchy did not develop his notion of the continuity of
a function before Bolzano developed his in RB, and that both notions
are essentially the same.  We argue that both assumptions are
incorrect, and that it is implausible that Cauchy's initial insight
into that notion, which eventually evolved to an approach using
infinitesimals, could have been borrowed from Bolzano's work.
Furthermore, we account for Bolzano's interest in that notion and
focus on his discussion of a definition by K\"astner (in Section 183
of his 1766 book), which the former seems to have misrepresented at
least partially.

Cauchy's treatment of continuity goes back at least to his 1817 course
summaries, refuting a key component of Grattan-Guinness' plagiarism
hypothesis (that Cauchy may have lifted continuity from RB after
reading it in a Paris library in 1818).  We explore antecedents of
Cauchy and Bolzano continuity in the writings of K\"astner and earlier
authors.

Keywords: Bolzano; Cauchy; K\"astner; continuity; infinitesimals;
variables 01A55; 26A15
\end{abstract}

\maketitle

\tableofcontents

\section{Introduction}
\label{sec1}

The issue of priority for the definition of the continuity of a
function was raised in \cite{Gr70} in a way that provoked controversy.
With regard to this issue, Grabiner seeks to shift the focus of
attention away from the Bolzano/Cauchy priority debate, and broaden
the discussion to include an analysis of their common predecessors,
particularly Lagrange.  She detects an ``immediate source of the
independent Bolzano--Cauchy definitions" both in Lagrange's 1798 book
\emph{Trait\'e de la r\'esolution des \'equations num\'eriques de tous
les degr\'es} and in his \emph{Th\'eorie des fonctions analytiques}
(see \cite[p.\;113]{Gr84}).  Grabiner concludes that ``these two books
are the most likely sources for both Cauchy's and Bolzano's
definitions of continuous function" (op.\;cit., p.\;114).  Grabiner's
analysis challenges Grattan-Guinness' claim that ``[Bolzano's and
Cauchy's] new foundations, based on limit avoidance, certainly swept
away the old foundations, founded largely on faith in the formal
techniques'' \cite[p.\;382]{Gr70}.  For sources of Bolzano's notion of
continuity in Lagrange see also \cite[p.\;422]{Ru99}.

Schubring similarly rules out Grattan-Guinness' hypothesis, and
furthermore challenges a common assumption that Bolzano's work was
virtually unknown in the mathematical community during the first half
of the 19th century \cite{Sc93}.  He reports on a (formerly) unknown
review of Bolzano's three important papers from 1816 and 1817, written
by a mathematician named J.\;Hoffmann in 1821 and published in 1823.

As for the Bolzano--Cauchy continuity, Grattan-Guinness investigated
the possibility of its antecedents, focusing on the following three
sources: (1) Cauchy's work prior to 1821, (2) Legendre, and (3)
Fourier; see \cite[p.\;286]{Gr70}.  His search reportedly did not turn
up any reasonable antecedents: ``of the new ideas that were to achieve
that aim -- of them, to my great surprise, I could find nothing''
(ibid.).  His investigation led him to his well-known controversial
conclusions.  What he missed were the following sources: (1) Cauchy's
earlier course summaries that were only discovered over a decade after
Grattan-Guinness' article (see Section~\ref{s23}); (2) Lagrange (as
argued by Grabiner); and (3) other 18th century authors, such as
K\"astner and Karsten (see Section~\ref{s4b}).

Some mathematicians and historians of mathematics assume that
Bolzano's definition of the continuity of a function in his 1817
\emph{Rein analytischer Beweis} preceded Cauchy's, and that the latter
first gave one in his 1821 textbook \emph{Cours d'Analyse}.  Both
assumptions turn out to be incorrect.  Scholars commonly assume the
following claims to be true:
\begin{enumerate}
[label={(Cl\,\theenumi})]
\item
\label{i1}
Bolzano and Cauchy gave essentially the same definition of
continuity, and 
\item
\label{i2}
Bolzano gave it earlier.
\end{enumerate}
We give some examples below.
\begin{itemize}
\item
Jarn\'{\i}\hspace{.5pt}k: ``Bolzano defines continuity essentially in
the same way as Cauchy does a little later'' \cite[p.\;36]{Ja81}.
\item
Segre: ``This led [Bolzano], in his \emph{Rein analytischer Beweis}
(written in 1817, four years before Cauchy published his \emph{Cours
d'analyse}), to give a definition of continuity and derivative very
similar to Cauchy's, etc.''  \cite[p.\;236]{Se94}.
\item
Ewald: ``[Bolzano's] definition is essentially the same as that given
by Cauchy in his \emph{Cours d'analyse} in 1821; whether Cauchy knew
of Bolzano's work is uncertain'' \cite[p.\;226]{Ew96}.
\item
Heuser: ``Cauchy defines continuity substantially in the same way as
Bolzano: \ldots''%
\footnote{In the original German: ``Stetigkeit definiert Cauchy
inhaltlich so wie Bolzano'' \cite[p.\;691]{He02}.  Heuser goes on to
present Cauchy's first 1821 definition in terms of $f(x+\alpha)-f(x)$
(see Section~\ref{s22b}), but fails to mention the fact that Cauchy
describes $\alpha$ as an \emph{infinitely small increment}.}
\end{itemize}
Now claim\;\ref{i1} is problematic since, as noted by L\"utzen,
\begin{quote}
Bolzano did not use infinitesimals%
\footnote{Note, however, that Bolzano did exploit infinitesimals in
his later writings; see e.g., \cite[note\;29, p.\;379]{Gr70},
\cite{Tr18}, and \cite{Fi20}.}
in his definition of continuity.  Cauchy did.  \cite[p.\;175]{Lu03}
\end{quote}
L\"utzen's claim that Cauchy used infinitesimals in his definition of
continuity is not entirely uncontroversial.  While Cauchy indisputably
used the term \emph{infiniment petit}, the meaning of Cauchy's term is
subject to debate.  Judith Grabiner \cite{Gr81}, Jeremy Gray
\cite[p.\;36]{Gr15}, and some other historians feel that a Cauchyan
infinitesimal is a sequence tending to zero.  Others argue that there
is a difference between null sequences and infinitesimals in Cauchy
(see e.g., \cite{19a}).

In sum, Cauchy's 1821 definitions exploited infinitesimals (and/or
sequences), whereas Bolzano's definition in \emph{Rein analytischer
Beweis} exploited the clause ``provided~$\omega$ can be taken as small
as we please'' in a way that can be interpreted as an incipient form
of an~$\epsilon,\delta$ definition relying on implied alternations of
quantifiers.  Such manifest differences make it difficult to claim
that the definitions were ``essentially the same.''

To determine the status of claim\;\ref{i2}, we will examine the
primary sources in Bolzano and Cauchy and compare their dates.

\section{Evolution of Cauchy's ideas documented by Guitard}
\label{s23}

\begin{figure}
\begin{center}
\includegraphics [scale=0.53,rotate=0,trim = 0 0 0 0,clip]
{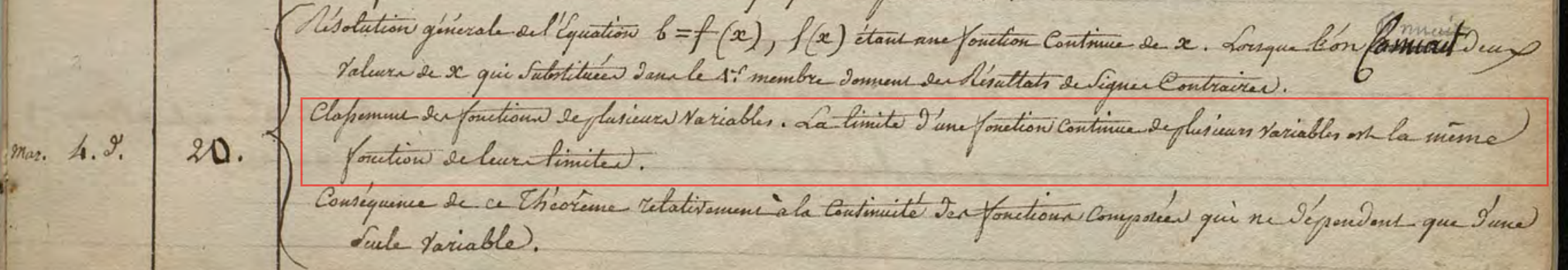}
\end{center}
\caption{Cauchy's treatment of continuity dating from 4\;march\;1817
in the gregorian calendar (which was a tuesday).  The ``Mar.''\;in the
figure stands for \emph{mardi}, tuesday.  The glyph resembling
$\partial$ to the right of the date seems to be shorthand for
\emph{ditto}, referring to the month of march mentioned on earlier
lines in this \emph{Registre de l'Instruction} for\;1817.}
\label{k19}
\end{figure}

Primary sources published in the 1980s suggest that an evolution took
place in Cauchy's ideas concerning continuity.  On 4\;march\;1817,
Cauchy presented an infinitesimal-free treatment of continuity in
terms of variables which is procedurally identical with the modern
definition of continuous functions via commutation of taking limit and
evaluating the function, as we discuss in Section~\ref{f3}.

\subsection{Continuity in 1817}
\label{f3}

In modern mathematics, a real function~$f$ is continuous at~$c\in\R$
if and only if for each sequence~$(x_n)$ converging to~$c$, one has
$f(\lim_{n\to\infty}x_n)=\lim_{n\to\infty}f(x_n)$, or briefly~$f\circ
\lim=\lim\circ f$ at~$c$.%
\footnote{The equivalence of such a definition with the
$\epsilon,\delta$ one requires the axiom of choice.}

In 1817, Cauchy wrote (see Figure~\ref{k19}):
\begin{quote} 
\emph{La limite d'une fonction continue de plusieurs variables est la
m\^eme fonction de leur limite}.  Cons\'equence de ce Th\'eor\`eme
relativement \`a la continuit\'e des fonctions compos\'ees qui ne
d\'ependent que d'une seule variable.%
\footnote{Translation: ``The limit of a continuous function of several
variables is [equal to] the same function of their limit.
Consequences of this Theorem with regard to the continuity of
composite functions dependent on a single variable.''  The reference
for this particular lesson in the Archives of the Ecole Polytechnique
is as follows: Le\;4\;Mars 1817, la le\c con 20.  Archives E. P., X II
C7, Registre d'instruction 1816--1817.}
(Cauchy as quoted in \cite[p.\;34]{Gu86}; emphasis added;
cf.\;\cite[p.\;255, note\;6 and p.\;309]{Be91})
\end{quote}
The Intermediate Value Theorem is proved in the same lecture.
Cauchy's treatment of continuity in 1817%
\footnote{Belhoste places it even earlier, in 1816: ``according to
the \emph{Registres}, Cauchy knew the modern concept of continuity as
far back as March 1817, but the `invention' was anterior, as shown by
the instructional program of December 1816'' \cite[p.\;255,
note\;6]{Be91}.}
contrasts with his definitions based on infinitesimals given four
years later in \emph{Cours d'Analyse} (CdA).  

\subsection{Continuity in \emph{Cours d'Analyse}}  
\label{s22b}

In CdA, Cauchy defines continuity as follows (see Figure~\ref{p1}):
\begin{figure}
\begin{center}
\includegraphics
[scale=0.25,rotate=0,trim = 0 0 0 0,clip]
{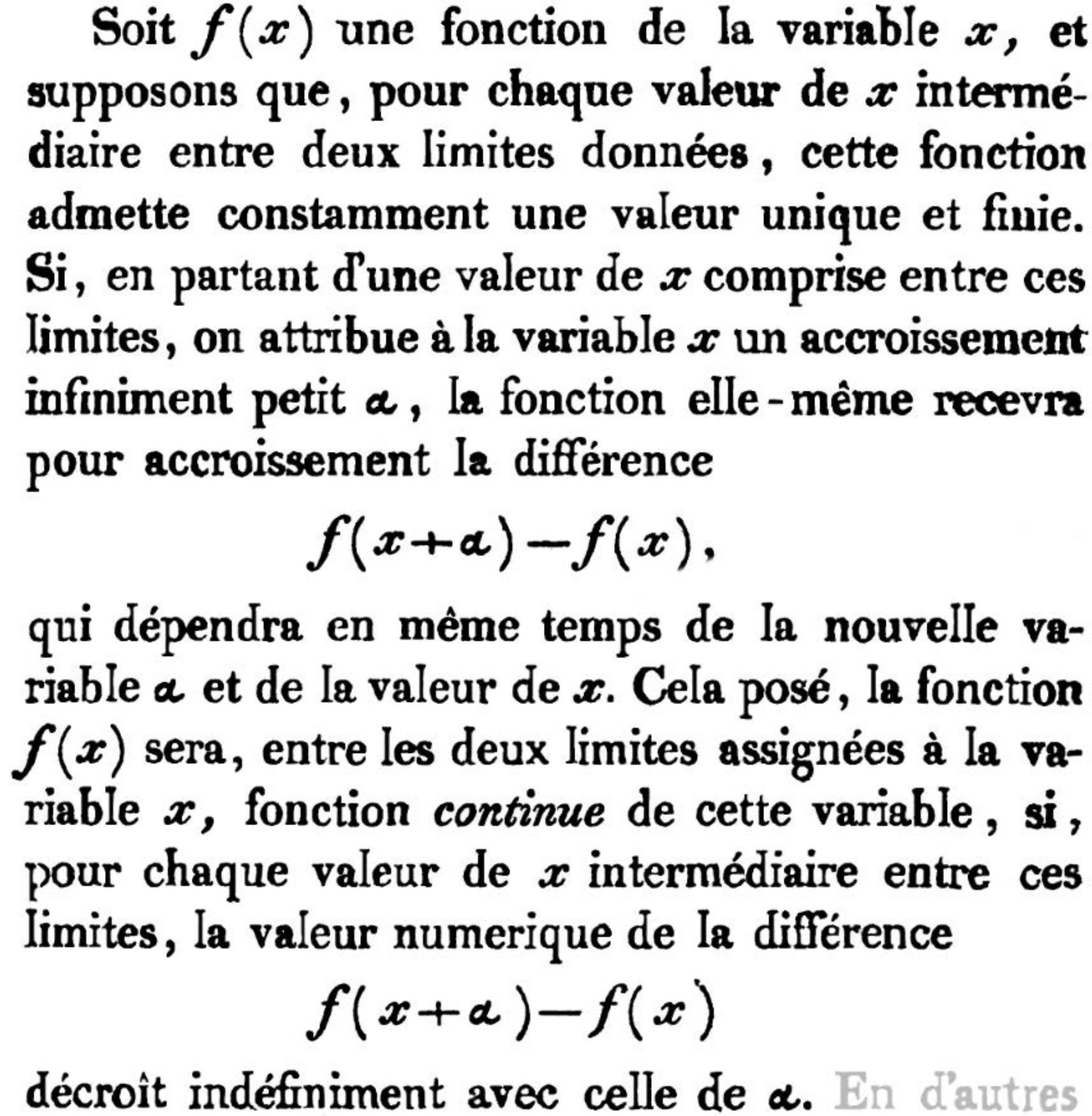}    

\end{center}
\caption{Cauchy's first 1821 definition of continuity}
\label{p1}
\end{figure}
\begin{quote}
Among the objects related to the study of infinitely small quantities,
we ought to include ideas about the continuity and the discontinuity
of functions.  In view of this, let us first consider functions of a
single variable.  Let~$f(x)$ be a function of the variable~$x$, and
suppose that for each value of~$x$ between two given limits, the
function always takes a unique finite value.  If, beginning with a
value of~$x$ contained between these limits, we add to the
variable~$x$ an \emph{infinitely small increment}~$\alpha$, the
function itself is incremented by the difference
\hbox{$f(x+\alpha)-f(x)$,} which depends both on the new
variable~$\alpha$ and on the value of~$x$.  Given this, the
function~$f(x)$ is a \emph{continuous} function of~$x$ between the
assigned limits if, for each value of~$x$ between these limits, the
numerical value of the difference~$f(x+\alpha)-f(x)$ decreases
indefinitely with the numerical value of~$\alpha$.  (Cauchy as
translated in \cite[p.\;26]{BS};%
\footnote{Reinhard Siegmund-Schultze writes: ``By and large, with few
exceptions to be noted below, the translation is fine'' \cite{Si}.}
emphasis on ``continuous'' in the original;
emphasis on ``infinitely small increment'' added)
\end{quote}
This definition can be thought of as an intermediary one between the
march\;1817 definition purely in terms of variables and containing no
mention of the infinitely small, and his second 1821 definition stated
purely in terms of the infinitely small (see Section~\ref{s22}).

\subsection{Second definition of continuity in CdA}
\label{s22}

\begin{figure}
\begin{center}
\includegraphics
[scale=0.25,rotate=0,trim = 0 0 0 0,clip]
{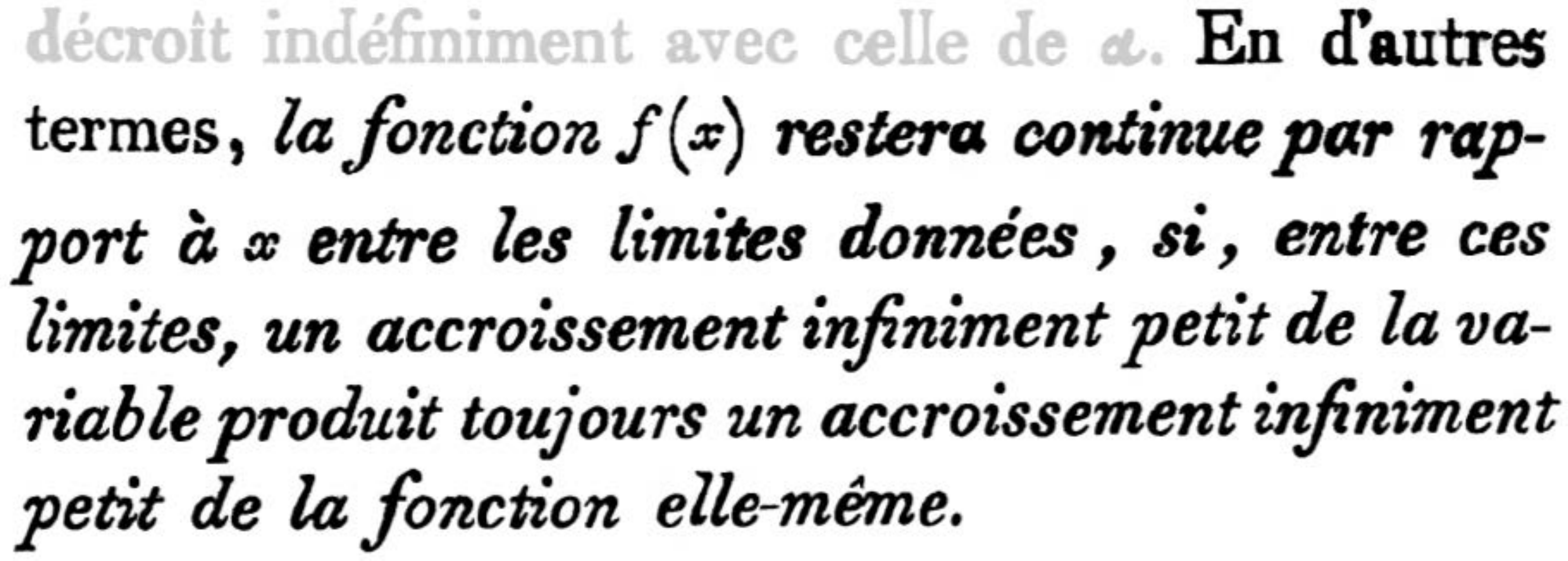}   
\end{center}
\caption{Cauchy's second 1821 definition of continuity}
\label{p2}
\end{figure}

Cauchy goes on to summarize the definition given above as follows (see
Figure~\ref{p2}):
\begin{quote}
In other words, \emph{the function~$f(x)$ is continuous with respect
to~$x$ between the given limits if, between these limits, an
infinitely small increment in the variable always produces an
infinitely small increment in the function itself}.%
\footnote{In the original: ``En d'autres termes, \emph{la fonction
$f(x)$ restera continue par rapport \`a~$x$ entre les limites
donn\'ees, si, entre ces limites, un accroissement infiniment petit de
la variable produit toujours un accroissement infiniment petit de la
fonction elle-m\^eme}'' \cite[pp.\;34--35]{Ca21}.}
(ibid.; emphasis in the original).
\end{quote}
Since Cauchy prefaced his second definition with the words \emph{en
d'autres termes} (``in other words"), he appears to have viewed the
pair of 1821 definitions as being equivalent.  Cauchy sums up his
discussion of continuity in CdA as follows:
\begin{quote}
We also say that the function~$f(x)$ is a continuous function of the
variable~$x$ in a neighborhood of a particular value of the
variable~$x$ whenever it is continuous between two limits of~$x$ that
enclose that particular value, even if they are very close together.
Finally, whenever the function~$f(x)$ ceases to be continuous in the
neighborhood of a particular value of~$x$, we say that it becomes
discontinuous, and that there is \emph{solution%
\footnote{meaning \emph{dissolution}, i.e., absence (of continuity).}
of continuity} for this particular value.  (ibid.; emphasis in the
original)
\end{quote}
Note that none of the 1821 definitions exploited the notion of limit.
We therefore find it puzzling to discover the contrary claim in a
recent historical collection:
\begin{quote}
Cauchy gave a faultless definition of continuous function, using the
notion of `limit' for the first time.  Following Cauchy's idea,
Weierstrass popularized the~$\epsilon$-$\delta$ argument in the 1870s,
etc.  \cite[p.\;283]{Da19}
\end{quote}
In a related vein, V\"ath opines that ``formulat[ing] properties which
hold for infinitesimals (which have been used by Leibniz) in an
$\epsilon$-$\delta$-type manner {\ldots} was first propagated by
Cauchy'' \cite[p.\;74]{Va07}.  Similarly, Goldbring and Walsh claim
the following:
\begin{quote}
[T]he mathematical status of [infinitesimals] was viewed as suspect
and the entirety of calculus was put on firm foundations in the
nineteenth century by the likes of Cauchy and Weierstrass, to name a
few of the more significant figures in this well-studied part of the
history of mathematics.  The innovations of their
``$\epsilon$-$\delta$ method'' {\ldots} allowed one to give rigor to
the na\"{\i}ve arguments of their predecessors.  \cite[p.\;843]{Go19}
\end{quote}
Presentist views of this type are, alas, not the exception, and much
work is required to counter them.  Recent work on Cauchy's stance on
the infinitely small and their applications includes \cite{17a},
\cite{17d}, \cite{18e}, and \cite{20a}.

To summarize, in 1817 Cauchy gave a characterisation of continuity in
terms of variables, whereas the second 1821 definition involved only
infinitesimals.  Meanwhile, the first 1821 definition exploited both
variables and infinitesimals.

\section
{Bolzano's \emph{Rein analytischer Beweis}}

Could Bolzano's \emph{Rein analytischer Beweis} (RB) \cite{Bo17} have
influenced Cauchy's definition of continuity?  Grattan-Guinness
wrote: 
\begin{quote}
Bolzano had given his paper [RB] two opportunities for publication,
for not only did he issue it as a pamphlet in 1817, but -- with the
same printing -- inserted it into the 1818 volume of the Prague
Academy \emph{Abhandlungen}.  That journal was available in Paris:
indeed, the \emph{Biblioth\`eque Imp\'eriale} (now the
\emph{Biblioth\`eque Nationale}) began to take it with \emph{precisely
the volume containing Bolzano's pamphlet}.  \cite[p.\;396]{Gr70}
(emphasis in the original)
\end{quote}
Of particular interest to us is Grattan-Guinness' reliance on the
availability of RB in the Paris \emph{Imperial Library} in 1818; see
Section~\ref{s41b}.  The papers \cite{Fr71} and \cite{Si73} provided
evidence against Grattan-Guinness' hypothesis.  However, as noted by
Jan Sebes\-tik, their work does not rule out the possibility that
``Cauchy could have read Bolzano's \emph{Rein analytischer Beweis} (or
heard about it) and could have been inspired by it'' \cite[pp.\;109,
111]{Se92}.  Thirty years after the Benis-Sinaceur paper, Russ wrote:
\begin{quote}
There has been \emph{discussion} in the literature on the possibility
that Cauchy might have plagiarized from Bolzano. See Grattan-Guinness
(1970), Freudenthal (1971) and Sinaceur (1973).  (\cite [p.\;149]
{Ru04}; emphasis added)
\end{quote}
It is our understanding that referring to the issue as a
``discussion'' tends to imply that the hypothesis of plagiarism has
not been definitively refuted.%
\footnote{Similarly, in a recently published book, Rusnock and
\v{S}ebest\'{\i}k mention that ``there has been speculation that
Cauchy may have learned a thing or two from Bolzano''
\cite[p.\;49]{Ru19}; see also note\;3 there.}
Arguably, therefore, the issue continues to have relevance.

\subsection{Grattan-Guinness' hypothesis}
\label{s41b}

Having summarized the historical background, Grattan-Guinness proceeds
to state his hypothesis:
\begin{quote}
So here is at least one plausible possibility for Cauchy to have found
a copy of Bolzano's paper, quite apart from the book-trade: he could
have noticed a new journal in the library's stock and examined it as a
possible course%
\footnote{Grattan-Guinness apparently means ``source.''}
of interesting research.  \cite[p.\;396]{Gr70}
\end{quote}
Grattan-Guinness specifically includes the concept of continuity in
his hypothesis (op.\,cit., p.\;374).  

It is our understanding that, while the evidence provided in the
articles \cite{Fr71} and \cite{Si73} shows clear and profound
differences between Cauchy and Bolzano's stance, it does not entirely
refute the aforementioned hypothesis.  We will provide a refutation of
a key component of Grattan-Guinness' hypothesis concerning the concept
of continuity.  Our refutation is based on the facts of the chronology
of the relevant works.  Namely, we will show that Cauchy possessed a
concept of continuity
\begin{enumerate}
\item
earlier than the date of the acquisition of a journal version of RB by
the \emph{Imperial Library} in Paris, and 
\item
even earlier than, or at least contemporaneously with, the date of the
Leipzig fair where RB was first marketed.
\end{enumerate}

Note that, according to Grattan-Guinness, the \emph{Biblioth\`eque
Imp\'eriale} started to take the journal where RB appeared in the year
1818.  Reading the 1818 journal version of RB could not therefore have
influenced Cauchy's treatment of continuity in 1817%
\footnote{\label{f1}Cauchy had discussed continuity even earlier, in
an 1814 article on complex functions (see \cite[p.\;380]{Fr71}).
However, that discussion stayed at the intuitive level and cannot be
described as reasonably precise.}
(see Section~\ref{s23}).  This refutes a key component of the
plagiarism hypothesis as proposed in \cite{Gr70} with regard to the
concept of continuity.  The comparison of dates establishes that
Cauchy's initial insight into continuity could not have been borrowed
from Bolzano's RB, though it does not rule out the possibility that
Cauchy may have been acquainted with Bolzano's work before formulating
the later, 1821 definitions in CdA.

Grattan-Guinness also brought broader plagiarism charges against
Cauchy, which are not refuted by our comparison of dates.  Notice,
however, that it is implausible that Cauchy may have seen Bolzano's
1816 text \emph{Der binomische Lehrsatz} \cite{Bo16}, where the latter
also gave a definition of continuity, since there is no evidence that
this text was available in France.  It seems that this is why
Grattan-Guinness found it necessary to speculate specifically
concerning the version of Bolzano's RB available in a Paris library in
1818, so as to bolster the plausibility of the plagiarism claim.
Grattan-Guinness may have had more of a point with regard to
E.\;G.\;Bj\"orling.  Apparently in the 1850s, Cauchy may not have been
transparent about possible influence of Bj\"orling's ideas related to
uniform convergence.  The issue was studied in \cite{Br07}.  For an
analysis of Cauchy's 1853 approach to uniform convergence see
\cite{18e}.

\subsection
{Bolzano's definition in \emph{Rein analytischer Beweis}}
\label{s32}

In his RB, Bolzano criticized some proofs of IVT for polynomials that
from his stance were ``based on an incorrect concept of
\emph{continuity},'' given for example their use of ``a truth borrowed
from \emph{geometry}'' or ``the introduction of the concepts of
\emph{time} and \emph{motion} \cite[pp.\;6, 8--9, 11]{Bo17}.  Instead,
he defined continuity as follows:
\begin{quote}
According to a correct definition, the expression that a function~$fx$
varies according to the law of continuity for all values of~$x$ inside
or outside certain limits means only that, if~$x$ is any such value
the difference~$f(x+\omega)-fx$ can be made smaller than any given
quantity, provided~$\omega$ can be taken as small as we please or (in
the notation we introduced in \S14 of \emph{Der binomische Lehrsatz}
etc., Prague, 1816)\;~$f(x+\omega)=fx+\Omega$.  (Bolzano as translated
in \cite[p.\;149, 256]{Ru04})
\end{quote}

The dating of RB will be analyzed in Section~\ref{s4} below.
Bolzano's definition is reasonably precise, as is Cauchy's approach.
Here ``reasonably precise'' means ``easily transcribable as a modern
definition'' (rather than merely an intuitive notion of continuity).%
\footnote{Note that we take no position with regard to which
definition was closer to a modern one, Bolzano's or Cauchy's
(Bolzano's was arguably closer to the modern \emph{Epsilontik}
standard).  The point we are arguing is that both were reasonably
precise in the sense specified.}
A modern formalisation of Bolzano's 1817 definition would involve
alternating quantifiers, whereas a modern formalisation of Cauchy's
1817 definition would retain almost verbatim the commutation of (a)
evaluating~$f$ and\;(b)\;taking~$\lim$ (see Section~\ref{f3}).
Apparently neither Jarn\'{\i}\hspace{.5pt}k nor Ewald (see
Section~\ref{sec1}) were aware of Cauchy's treatment of both
continuity and the IVT dating from 4\;march\;1817.

\subsection{The dating of Bolzano's RB}
\label{s4}

The earliest known written record of Bolzano's RB is in a catalog of
the Easter book fair at Leipzig. 

\begin{figure}
\begin{center}
\includegraphics
[scale=0.25,rotate=0,trim = 0 0 0 0,clip]
{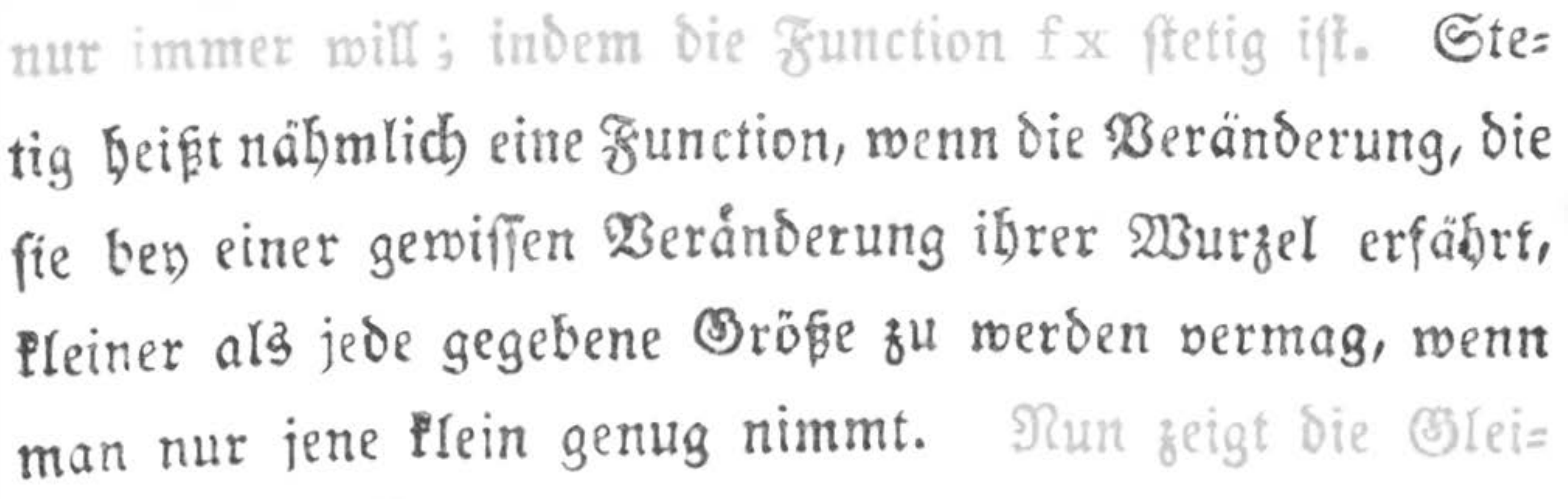}     
\end{center}
\caption{Bolzano's definition of continuity}
\label{p3}
\end{figure}

According to \cite[p.\;4]{Ev14}, both the catalog \cite[p.\;30]{Ol17}
and the fair itself date from 27\;april\;1817, over a month later than
the earliest written record of Cauchy's treatment of continuity.  It
should be noted, however, that Bolzano also gave a definition of
continuity in an 1816 publication \cite{Bo16} (see Figure~\ref{p3}):
\begin{quote}
For a function is called continuous if the change which occurs for a
certain change in its argument, can become smaller than any given
quantity, provided that the change in the argument is taken small
enough.%
\footnote{In the original: ``Stetig hei{\ss}t n{\"a}hmlich eine
Function, wenn die Ver{\"a}nderung, die sie bey einer gewissen
Ver{\"a}nderung ihrer Wurzel erf{\"a}hrt, kleiner als jede gegebene
Gr{\"o}{\ss}e zu werden vermag, wenn man nur jene klein genug nimmt''
\cite[p.\;34]{Bo16}.  Note that Bolzano repeatedly uses \emph{Wurzel}
in the sense of ``input to a function''; see e.g., footnote on
page\;11 of \cite{Bo17}.  The issue is discussed in \cite[p.\;256,
note f]{Ru04}.}
(Bolzano as translated in \cite[p.\;184]{Ru04})
\end{quote}

This definition is immediately followed by an attempted proof of an
erroneous assertion.  Namely, Bolzano claims to prove that if a
function~$F$ is differentiable then its derivative,~$f$, is
continuous.  This indicates that Bolzano's definition of continuity
was still sufficiently ambiguous to accomodate errors, as was
his~$\omega/\Omega$ notation.  Recently \cite[Abstract]{Fu19} have
argued in \emph{Historia Mathematica} that ``those quantities [i.e.,
Bolzano's~$\omega$] are not clearly `proto-Weierstrassian'.''

It is worth noting that an even earlier mention of ideas in the
direction of Bolzano's definition of continuity occurs in Bolzano's
mathematical diaries of early 1815: ``if therefore~$\xi$ is taken
smaller than any given quantity, i.e.~$=\omega$, the value
of~$f(x+\omega)-fx$ must be able to become as small as desired'' (see
op.\;cit., note\;86).  Insofar as Cauchy had no access either to
Bolzano's diaries or the latter's 1816 work, and the former would have
formulated his first definition of continuity shortly before or in any
case at about the same time as the 1817 Easter book fair at Leipzig,
it is implausible that Cauchy's 1817 definition could have been
borrowed from Bolzano's work.

\section{Antecedents in K\"astner, Karsten, and others}
\label{s4b}

There exists a historiographic controversy with regard to the issue of
continuity in the historical development of mathematics.  Unguru and
his disciples adopt a radical posture against such continuity.  Other
scholars endorse continuity at various levels and to varying extent.
We adopt the latter view, to the extent that we detect continuity
between, for example, the work of K\"astner, on the one hand, and that
of Bolzano and Cauchy, on the other.  For more details see \cite{20c}.

The mathematical diaries of Bolzano written during 1814--1815 also
contain criticism of, e.g., \cite{Ca97} and \cite{Cr13} because of
their assumption of the law of continuity: in the first case he stated
that in such a law ``[lay] the key for the resolution of the whole
riddle of infinitesimal calculus" \cite[p.\;152]{Bo95}; in the latter
case, he pointed out that K\"astner had ``already drawn attention to
the surreptitious acceptance of this law" \cite[p.\;144]{Bo97}.  As we
already mentioned, the first published record of a definition of
continuity given by Bolzano dates from the following year, after which
he published his reasonably precise definition included in RB.

\begin{figure}
\begin{center}
\includegraphics [scale=0.25,rotate=0,trim = 0 0 0 0,clip]
{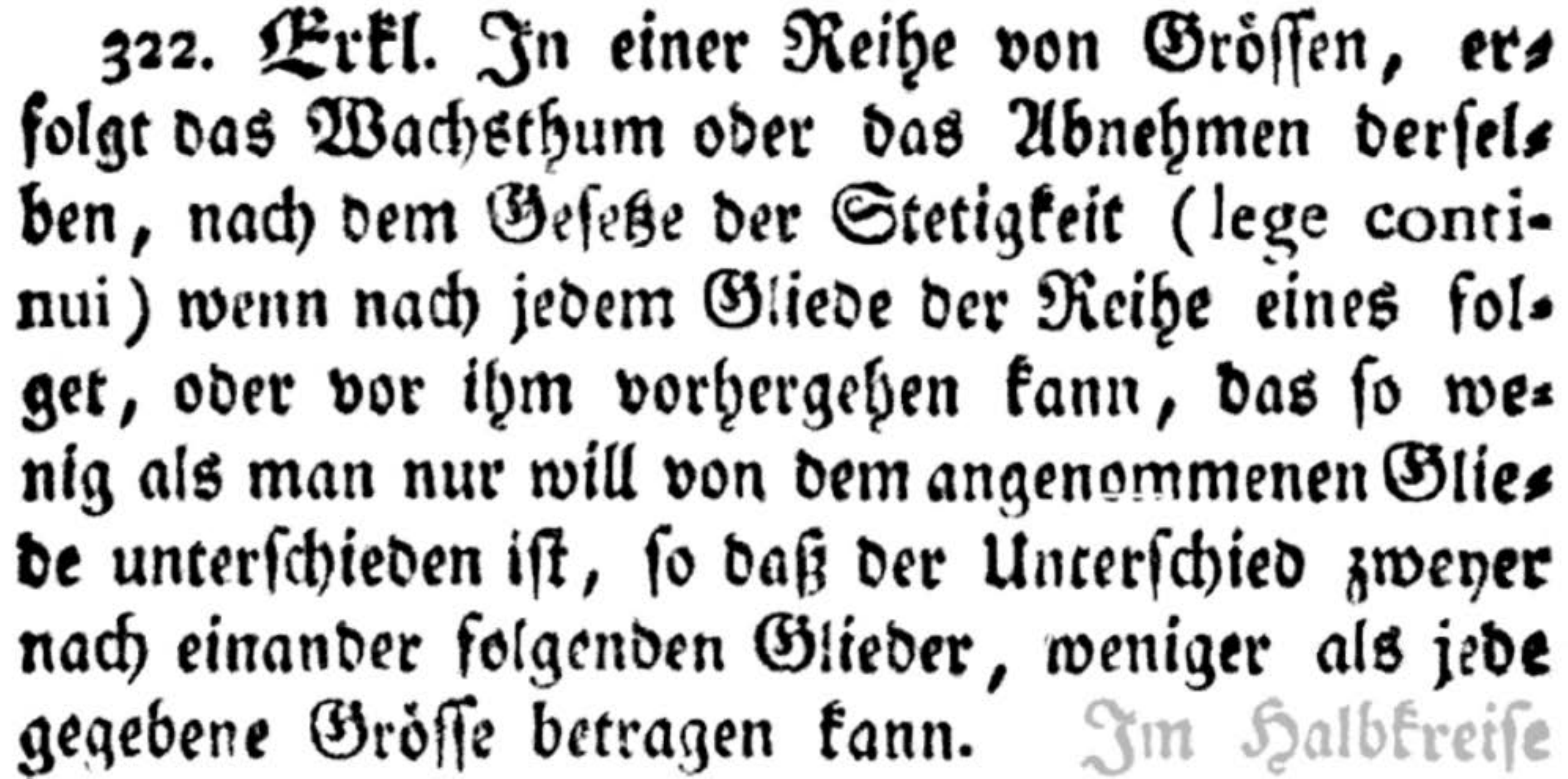}    
\end{center}
\caption{K\"astner's 1760 definition of continuity}
\label{p4}
\end{figure}

As his later works and mathematical diaries show, Bolzano continued to
be interested in that issue.  Thus, in his \emph{Theory of Functions},
written in the 1830s, he would have ``sharpened'' his 1817 definition
\cite[p.\;306]{Ru05}.  Rusnock and Kerr-Lawson argue that, as early as
the 1830s, Bolzano not only grasped the distinction between pointwise
continuity and uniform continuity but also presented a pair of key
theorems concerning the latter (ibid.).  Moreover, in that work
Bolzano acknowledged that ``[t]he concept of continuity has already
been defined essentially as I do here by [other contemporary authors]"
such as Cauchy and Ohm \cite[p.\;449]{Ru04}.  However, at the same
time, in that work he criticized certain specific definitions,
including one by A. G. K\"astner in 1766.  On the one hand, Bolzano's
definition surely constitutes an improvement upon the definition of
local continuity by K\"astner in 1760 (see Figure~\ref{p4}).  On the
other hand, Bolzano seems to have misrepresented, at least partially,
the relevant passage from K\"astner's work of 1766.

\subsection{K\"astner's 1760 definition}
\label{s41}

K\"astner's definition included in his volume on the analysis of
finite quantities (\emph{Analysis endlicher Gr\"ossen}), or
letter-algebra, and which can be found in a section entitled ``On
curved lines,'' runs as follows:
\begin{quote}
In a sequence%
\footnote{We translated \emph{Reihe} as `sequence', even though it is
often translated as `series', since `series' nowadays is a standard
technical term which is not appropriate here, and moreover the German
term \emph{Reihe} can mean either `sequence' or `series'.}
of magnitudes, their increase or decrease takes place in accordance
with the law of continuity (lege continui), if after each term of the
sequence, another one follows or precedes \emph{the given term} that
differs from it [i.e., from the given term] by as little as one
wishes; as a consequence,%
\footnote{The German conjunction \emph{so dass}, especially in
K\"astner's (now obsolete) spelling as two separate words, resembles
the English `such that'; in the present case, however, this is a false
friend.  In fact `as a consequence' is one of several standard
translations of the German conjunction \emph{sodass}.}
the difference of two consecutive terms%
\footnote{K{\"a}stner's phrasing \emph{nach einander folgender} could
possibly be interpreted as the statement that the terms mentioned here
are immediate successor elements, in particular since the standard
technical translation for `immediate sucessor element' is
\emph{Nachfolger}.  This, however, could not be what K{\"a}stner meant
to say.  K\"astner's phrasing (note that he does not say
\emph{Nachfolger} outright) is sufficiently vague to allow for an
interpretation where he means to speak of two terms which follow
\emph{shortly} one after another, though there are other terms in
between.}
can amount to less than any given magnitude.%
\footnote{In the original: ``In einer Reihe von Gr{\"o}ssen, erfolgt
das Wachsthum oder das Abnehmen derselben, nach dem Gesetze der
Stetigkeit (lege continui) wenn nach jedem Gliede der Reihe eines
folget, oder vor ihm vorhergehen kann, das so wenig als man nur will
von dem angenommenen Gliede unterschieden ist, so da{\ss} der
Unterschied zweyer nach einander folgender Glieder, weniger als jede
gegebene Gr{\"o}sse betragen kann.''  This was quoted in
\cite[p.\;283]{Sp15}.  In our translation, we tried to strike a
balance between literalness and readability in line with an approach
taken in \cite{Bl18}.}
\cite[paragraph\;322, p.\;180]{Ka60}
\end{quote}

\subsection{K\"astner's influence on Bolzano}

Russ notes K\"astner's influence on Bolzano in the following terms:
\begin{quote}
[T]here were two authors, Wolff and K\"astner, whose work, between
them, dominated the century in the Ger\-man-speaking regions.
\ldots~[T]hey were both committed to education and wrote highly
systematic and comprehensive multivolume textbooks on mathematics that
went through many editions and were very influential.  Not
surprisingly, they were both authors to whom Bol\-zano makes frequent
reference in his early works. \\ \cite[p.\;14]{Ru04}
\end{quote}
Indeed, in Bolzano's mathematical diaries there is a note from the
early 1820s, entitled ``On the law of continuity.''  Bolzano's note
includes a reference to paragraph 183 of K\"astner's work on mechanics
\cite{Ka66} and to paragraph 235 of W. J. G. Karsten's work on
mechanics \cite{Ka69}; see \cite[p.\;63]{Bo05}.  The formulation of
both authors ultimately relied on the notion of continuity according
to which ``[a] \emph{continuous quantity} (continuum) is that
[quantity] whose parts are all connected together in such a way that
where one ceases, another immediately begins, and between the end of
one and the beginning of another there is nothing that does not belong
to this quantity'' \cite[p.\;17]{Ru04}; see \cite[p.\;209]{Ka67}; but
only that of Karsten would be equivalent to IVT \cite[p.\;223]{Ka69}.
Interestingly, as we already mentioned, in a later work Bolzano went
back to discuss the notion of continuity in that paragraph of
K\"astner's work.  We will analyze such a reception of the latter's
ideas in Section~\ref{s43}.

\subsection
{Bolzano misattributes a definition to K\"astner}
\label{s43}

We reviewed K\"astner's 1760 definition in Section~\ref{s41}.  In his
\emph{Theory of Functions}, Bolzano seems to have mistakenly
attributed a different definition to K\"astner in 1766, which he
(Bolzano) considered to be ``too broad'':
\begin{quote}
Some very respected mathematicians like \emph{K\"astner}
(\emph{h\"ohere Mechanik}, Auflage 2, \S\S\,183 ff.)  and \emph{Fries}
(\emph{Naturphilosophie}, \S\,50) define the continuity of a function
$Fx$ as that property of it by virtue of which it does not go from a
certain value~$Fa$, to another value~$Fb$, without first \emph{having
taken all the values lying in between}.  However, it will be seen
subsequently that this definition is \emph{too wide}%
\footnote{Perhaps a better translation is ``too broad''.}
if in fact the concept intended is to be equivalent to the one above.
(Bolzano as translated in \cite[p.\;449]{Ru04}; emphasis on K\"astner
and Fries in the original; emphasis on ``having taken all the values
lying in between'' and ``too wide'' added)
\end{quote}
As we already noted, K\"astner's formulation to which Bolzano refers
here ultimately relied on the former's geometric notion of continuity.
So, while K\"astner's paragraph\;183 is part of a section ``On the law
of continuity'' (which in turn is part of a chapter ``On the movement
of solid bodies with determined magnitude and shape''), he explicitly
refers to the note in his definition 6 (straight and curved lines) of
his book on geometry.  In that note K\"astner points out that before
the curved line that goes from $A$ to $B$ reaches $B$, ``all the minor
changes in between must occur" \cite[p.\;161]{Ka58}.

Bolzano would seem to attribute a different definition (via the
satisfaction of the Intermediate Value Theorem) to K\"astner (as well
as to Fries) in the particular case of that paragraph.  Nonetheless,
Bolzano's attribution appears to be incorrect.

In fact, K\"astner's discussion of the law of continuity in his
section\;183 resembles, to some extent, Cauchy's definition of
continuity based on infinitesimals given in Section~\ref{s22} above
(though of course K\"astner's viewpoint is geometric rather than
analytic):
\begin{quote}
On the Law of Continuity.  183. In the investigation which we now
present, it is assumed that the speed of a body does not change
instantaneously, but rather by infinitely small gradations. Just the
same can be said of the direction.  If one views the matter from that
perspective, then a body which is being reflected does not change its
direction instantaneously to the opposite direction: its speed becomes
smaller and smaller in the previous direction, finally vanishes, and
then transforms into a velocity having the opposite direction. This is
the \textbf{Law of Continuity} (applied to these matters).  To wit, by
the latter law one claims that generally, no change happens suddenly,
but that \emph{every change always moves through infinitely small
gradations} (of which already the movement of a point along a curve is
an example; [cf. K\"astner's] Geom.\;6.\;Erkl.\;Anm.).
(\cite[p.\;350, \S\,183]{Ka66};%
\footnote{According to \cite[Abbildung 10]{Kr14}, there were two
edititions of this treatise.  These are \cite{Ka66} and \cite{Ka93}.
In the 1793 edition of K\"astner's treatise referred to by Bolzano as
\emph{Auflage\;2}, Section\;183 appears on page 543.}
emphasis on ``law of continuity'' on the original; emphasis on ``every
change, etc.''  added)
\end{quote}
What may have led Bolzano to claim that K\"astner defined continuity
based on the satisfaction of IVT?  Note that K\"astner's text contains
the following three sentences:
\begin{enumerate}
[label={(K\theenumi})]
\item
\label{s1}
If one views the matter from that perspective, then a body which is
being reflected does not change its direction instantaneously to the
opposite direction: its speed becomes smaller and smaller in the
previous direction, finally vanishes, and then transforms into a
velocity having the opposite direction.
\item
\label{s2}
This is the Law of Continuity (applied to these matters).
\item
\label{s3}
To wit, by the latter law one claims that generally, no change happens
suddenly, but that every change always moves through infinitely small
gradations.
\end{enumerate}
Possibly, Bolzano interpreted sentence~\ref{s1} as the definition of
the law of continuity mentioned in sentence~\ref{s2}.  Now
sentence~\ref{s1} does sound like (a physical interpretation of) the
IVT.

However, reading the three sentences together, it is clear that
K\"astner meant sentence~\ref{s3} to be the detailed formulation of
the law of continuity.  Meanwhile, in sentence~\ref{s2}, K\"astner
specifically uses the verb \emph{applied}.  This indicates that
K\"astner thinks of sentence~\ref{s1} as an \emph{application} of the
law of continuity, rather than the formulation thereof.  Now in modern
mathematics it is certainly true that continuity implies IVT, though
the converse is incorrect, as Bolzano himself argued (see
\cite[\S\,84, pp.\;471--472]{Ru04}).  In his \emph{Theory of
Functions}, Bolzano outlines an idea for a function that takes every
intermediate value without being continuous, as follows.

Bolzano starts with an everywhere discontinuous function $W(x)$
described in {\S}37, defined only on a collection of rational points,
and built out of a pair of linear functions of different slope.  In
{\S}39, Bolzano asserts that the remaining infinitely many points can
be used to assign the values of the function so as to ``fill in''
whatever values are missing.  Bolzano's argument is mentioned in
\cite[p.\;395]{Se92}%
\footnote{Sebestik also points out that Bolzano and Cauchy's
definitions of continuity could have been ``the result of a critical
reflection on the texts by Euler and Lagrange'' \cite[pp.\;110,
81--83]{Se92}.}
and \cite{Sm17} (see p.\;52 and note\;49 there).  For a study of
counterexamples to the implication ``if~$f$ satisfies IVT then~$f$ is
continuous'' see \cite{Om14}, \cite{Ra16}, and \cite{De18}.

In conclusion, Bolzano may have interpreted sentence~\ref{s1} as the
formulation of continuity (rather than an application thereof).
Unlike Cauchy, Bolzano seems never to have formulated a definition of
continuity in terms of infinitesimals.  It is possible that
K\"astner's sentence~\ref{s3} made no sense to Bolzano, who was
therefore led to take sentence~\ref{s1} to be the formulation of
continuity.  Thus, while Fries may perhaps have given a different
definition of continuity via the satisfaction of IVT (as Bolzano
claimed), K\"astner apparently did not.

\subsection{Continuity in Leibniz}

An even earlier source for local continuity may have influenced
K\"astner and other 18th century authors.  Such a source is in
Leibniz's 1687 formulation of the principle of continuity:
\begin{quote}
When the difference between two instances in a given series or
\emph{that which is presupposed} can be diminished until it becomes
smaller than any given quantity whatever, the corresponding difference
in \emph{what is sought} or in their results must of necessity also be
diminished or become less than any given quantity whatever.  (Leibniz
as translated by Loemker in \cite[p.\;351]{Le89}; emphasis added)
\end{quote}
In modern terminology, Leibnizian ``what is sought" is the dependent
variable while ``that which is presupposed" is the independent
variable.  What Leibniz refers to as the principle of continuity%
\footnote{Not to be confused with his \emph{law of continuity}.  For a
detailed discussion see \cite{13f}, \cite{14c}, \cite{16a},
\cite{17b}, \cite{18a}.}
involves, in modern terminology, the condition that a convergent
sequence in the domain should get mapped to a convergent sequence in
the range.%
\footnote{In modern analysis, the sequence-condition is equivalent to
continuity for first-countable spaces.}

Cauchy's approach dating from 4\;march\;1817 is not the final word on
continuity, but it can be described as reasonably precise (in the
sense explained in Section~\ref{s32}).  This is unlike many intuitive
definitions given earlier%
\footnote{Including Cauchy's own definition in 1814, in an article on
complex functions quoted by Freudenthal; cf.\;note~\ref{f1}.}
that cannot be so formalized.  

Notice that Bolzano's definition is similarly reasonably precise but
also not without its problems.  Thus, the~$\Omega$ appearing there
seems to be defined as the difference~$f(x+\omega)-f(x)$, whereas the
corresponding~$\epsilon$ in the modern definition is
a~$\forall$-quantified variable entirely unrelated to~$f$.  It is
possible that this was also Bolzano's intention, but it must be
admitted that such an intention was only imperfectly expressed by
Bolzano's formula~$f(x+\omega)=fx+\Omega$ and accompanying comments;
see \cite{Fu19} for a fuller discussion.

\section{Conclusion}

We have re-examined the priority issue with regard to the concept of
continuity.  Course notes available at the Ecole Polytechnique
indicate that Cauchy had a reasonably precise concept of continuity of
a function earlier than is generally thought.  In particular Cauchy's
concept was earlier than, or at least contemporaneous with, the first
written record of Bolzano's 1817 work \emph{Rein analytischer Beweis}.

In 1970, Grattan-Guinness speculated that Cauchy may have read a
version of Bolzano's \emph{Rein analytischer Beweis} found in a Paris
library in 1818, and subsequently plagiarized some of Bolzano's
insights, including continuity, when writing the 1821 \emph{Cours
d'Analyse}.  Such a hypothesis is refuted by a written record of a
reasonably precise treatment of continuity by Cauchy dating from march
1817, and hence anterior to the Paris library acquisition, on which,
among other things, Grattan-Guinness based his hypothesis.

The proximity of the dates indicates an independence of Cauchy's and
Bolzano's scientific insight, and should contribute not only to end
speculations as to possible plagiarism (with regard to the notion of
continuity) on either side, but also to improve our understanding of
their respective developments of such a notion.

The prototypes of both Bolzano's and Cauchy's definitions of
continuity in formulations found in 18th century and early 19th
century works, such as those of K\"astner, are yet to be explored
fully.

\section*{Acknowledgments} 
We thank the anonymous referees for numerous suggestions that helped
improve the article.  We are grateful to Olivier Azzola, archivist at
the Ecole Polytechnique, for granting access to Cauchy's handwritten
course summaries reproduced in Figure~\ref{k19}.  E.\;Fuentes
Guill\'en was supported by the Postdoctoral Scholarship Program of the
Direcci\'on General de Asuntos del Personal Acad\'emico (DGAPA-UNAM).
V. Kanovei was supported by RFBR grant no.\;18-29-13037.

\end{document}